\documentclass[a4paper,12pt]{article}

\usepackage{a4,amssymb,amsfonts,verbatim,pstricks,pst-plot,xypic,times}
\usepackage{yfonts}

\usepackage[frenchb,english]{babel}
\usepackage{amssymb}
\usepackage{amsmath}

\usepackage{amssymb,latexsym}
\usepackage{mltex}
\usepackage{amsmath}
\usepackage{enumerate}
\usepackage{amsthm}

\newtheorem{thm}{Theorem}     
     
\newtheorem{lem}{Lemma}[section]
\newtheorem{prop}[lem]{Proposition}     
\newtheorem{rem}{Remark}  

\newtheorem{defn}[lem]{Definition}

\newcommand{\gn}{\nabla} 
 
\newcommand{\Ekt}{\mathbb{E}(\kappa,\tau)}
\newcommand{\Ekto}{\mathbb{E}(\kappa_1,\tau_1)}
\newcommand{\Ektt}{\mathbb{E}(\kappa_2,\tau_2)}

\newcommand{\nablab}{\overline{\nabla}}

\newcommand{\Ric}{\mathrm{Ric}}

\newcommand{\PSN}{P_{\Ss^1}N}

\newcommand{\ov}{\overline}
\newcommand{\PSM}{P_{\Ss^1}M}

\newcommand{\iid}{\mathrm{Id}\,}

\newcommand{\spin}{\mathrm{Spin}}

\newcommand{\ddet}{\mathrm{det}\,}

\newcommand{\eend}{\mathrm{End}\,}

\newcommand{\trace}{\mathrm{tr\,}}

\newcommand{\C}{\mathbb{C}}

\newcommand{\<}{\left\langle}       
\renewcommand{\>}{\right\rangle}

\newcommand{\Ss}{\mathbb{S}}
\newcommand{\HH}{\mathbb{H}}

\newcommand{\R}{\mathbb{R}}

\newcommand{\M}{\mathbb{M}}

\newcommand{\Chi}{\mathfrak{X}}
\newcommand{\pre}{\Re e}

\newcommand{\spinc}{\mathrm{Spin^c}}
\newcommand{\CP}{\mathbb{C}P^2}
\newcommand{\CH}{\mathbb{C}H^2}
\newcommand{\CM}{\mathbb{M}_{\C}^2(c)}

\newcommand{\beqt}{\begin{equation}}  \newcommand{\eeqt}{\end{equation}}
\newcommand{\bal}{\begin{align}}      \newcommand{\eal}{\end{align}}
\newcommand{\ba}{\begin{array}}      \newcommand{\ea}{\end{array}}
\newcommand{\bc}{\begin{center}}     \newcommand{\ec}{\end{center}}
\newcommand{\be}{\begin{enumerate}}  \newcommand{\ee}{\end{enumerate}}
\newcommand{\beq}{\begin{eqnarray}}  \newcommand{\eeq}{\end{eqnarray}}
\newcommand{\beQ}{\begin{eqnarray*}} \newcommand{\eeQ}{\end{eqnarray*}}
\newcommand{\bi}{\begin{itemize}}    \newcommand{\ei}{\end{itemize}}
\newcommand{\bt}{\begin{tabular}}    \newcommand{\et}{\end{tabular}}

\title{Hypersurfaces of Spin$^c$ Manifolds and Lawson Type Correspondence}
\author{Roger Nakad and Julien Roth}
\date{\today}

\begin{document}
\maketitle
\begin{abstract} 
\noindent Simply connected $3$-dimensional homogeneous manifolds $\Ekt$, with $4$-dimen\-sional isometry group, have a 
canonical $\spinc$ structure carrying parallel or Killing spinors. The restriction to any hypersurface of these parallel or
 Killing spinors allows to characterize isometric immersions of surfaces into $\Ekt$. 
As application, we get an elementary proof of a  Lawson type correspondence for constant mean curvature 
surfaces in $\Ekt$. Real hypersurfaces of the complex projective space and the complex hyperbolic space are also characterized via $\spinc$ spinors. 
\end{abstract} 

\section{Introduction}
It is well-known that a conformal immersion of a surface in $\R^3$ could be characterized by a spinor field $\varphi$ satisfying
\beqt\label{eqdir1}
D\varphi=H\varphi,
\eeqt
where $D$ is the Dirac operator and $H$ the mean curvature of the surface (see \cite{KS} for instance). In \cite{Fr}, Friedrich characterized surfaces in $\R^3$ in a geometrically invariant way. More precisely, consider an isometric immersion of a surface $(M^2,g)$ into $\R^3$. The restriction to $M$ of a parallel spinor of $\R^3$ satisfies, for all $X\in \Gamma(TM)$, the following relation
\beqt\label{eqrestr}
\nabla_X\varphi=-\frac{1}{2}IIX\bullet\varphi,
\eeqt
where $\nabla$ is the spinorial Levi-Civita connection of $M$, ``$\bullet$'' denotes the Clifford multiplication of $M$ and $II$ is the shape operator of the immersion. Hence, $\varphi$ is a solution  of the Dirac equation \eqref{eqdir1} with constant norm. Conversely, assume that  a Riemannian surface $(M^2,g)$ carries a
spinor field $\varphi$, satisfying
\begin{eqnarray}
 \nabla_X \varphi = -\frac 12 EX \bullet \varphi,
\label{EEE}
\end{eqnarray}
where $E$ is a given symmetric endomorphism on the tangent bundle. It is straightforward to see that $E = 2 \ell^\varphi.$ Here $\ell^\varphi$ is a field of symmetric endomorphisms associated with the field of quadratic forms, denoted also by $\ell^\varphi$, called the energy-momentum tensor which is given, on the complement set of zeroes of $\varphi$, by 
$$\ell^\varphi (X) = \Re \< X\bullet\nabla_X \varphi,\frac{\varphi}{\vert\varphi\vert^2}\>,$$
for any $X\in \Gamma(TM)$. Then, the existence of a pair $(\varphi, E)$  satisfying (\ref{EEE}) implies that the tensor $E= 2\ell^\varphi$ satisfies the Gauss and Codazzi equations and by Bonnet's theorem, there exists a local isometric immersion of $(M^2,g)$ into $\R^3$ with $E$ as shape operator. Friedrich's result was extended by Morel \cite{Mo} for surfaces of the sphere $\Ss^3$ and the hyperbolic space $\HH^3$.\\

Recently, the second author \cite{Rot}  gave a spinorial characterization of surfaces isometrically immersed into  $3$-dimensional homogeneous manifolds with $4$-dimensional isometry group. These manifolds, denoted by $\Ekt$ are Riemannian fibrations over a simply connected 2-dimensional manifold $\M^2(\kappa)$ with constant curvature $\kappa$ and bundle curvature $\tau$. This fibration can be represented by  a unit vector field $\xi$ tangent to the fibers.\\

The manifolds $\Ekt$ are $\spin$ having a special spinor field $\psi$. This spinor is constructed using real or imaginary Killing spinors on $\M^2(\kappa)$. If $\tau \neq 0$, the restriction  of $\psi$ to a surface gives rise to a spinor field $\varphi$ satisfying, for every vector field $X$,
\begin{eqnarray}
\nabla_X\varphi=-\frac{1}{2}IIX\bullet\varphi+
i\frac{\tau}{2} X\bullet\overline\varphi -i\frac{\alpha}{2} g(X, T) T\bullet\overline \varphi + i\frac{\alpha}{2} f g(X, T) \ov \varphi.
\label{varr}
\end{eqnarray}
Here $\alpha= 2\tau -\frac{\kappa}{2\tau}$,  $f$ is a real function and  $T$ is a vector field on $M$ such that $\xi = T +f \nu$ is the decomposition of $\xi$ into 
tangential and normal parts ($\nu$ is the normal vector field of the immersion). The spinor  $\overline \varphi$ is given by $\ov \varphi := \varphi^+ - \varphi^-$, where $\varphi =\varphi^+ + \varphi^- $ is the decomposition into positive and negative spinors.  Up to some additional geometric assumptions on $T$ anf $f$, the spinor $\varphi$ allows to characterize the immersion of the surface into $\Ekt$ \cite{Rot}. \\

In the present paper, we  consider $\spinc$ structures on $\Ekt$ instead of $\spin$ structures. The manifolds $\Ekt$ have a  canonical $\spinc$ structure carrying a natural spinor field, namely a real Killing spinor with Killing constant $\frac{\tau}{2}$. The restriction of this Killing spinor to $M$ gives rise to a special spinor satisfying
$$\nabla_X\varphi=-\frac{1}{2}IIX\bullet\varphi+i\frac{\tau}{2}X\bullet\overline{\varphi}.$$
This spinor, with a curvature condition on the auxiliary bundle, allows the characterization of the immersion of $M$ into $\Ekt$ without any additional geometric assumption on $f$ or $T$ (see Theorem \ref{thmEkt}). From this characterization, we get an elementary spinorial proof of a Lawson type correspondence for constant mean curvature surfaces in $\Ekt$ (see Theorem \ref{thmdan}).\\

The second advantage of using $\spinc$ structures in this context is when we consider hypersurfaces
 of $4$-dimensional manifolds. Indeed, any oriented $4$-dimensional K\"ahler manifold has a canonical $\spinc$
 structure with parallel spinors. In particular, the complex space forms $\CP$ and $\CH$. Then, using an analogue 
of Bonnet's Theorem for complex space forms, we prove a spinorial characterization of hypersurfaces of the 
complex projective space $\CP$ and of the complex hyperbolic space $\CH$. 
This work generalizes to the complex case the results of \cite{Mo} and \cite{LR}. Finally, we apply this characterization for Sasaki hypersurfaces.
\section{Preliminaries}
In this section we briefly introduce basic facts about $\spinc$ geometry of hypersurfaces (see \cite{LM, montiel, fr1, r1, r2}). Then we give a short description of the complex space form $\CM$ of complex dimension $2$, the $3$-dimensional homogeneous manifolds with $4$-dimensional isometry group $\Ekt$ and their hypersurfaces (see \cite{Da2, Sco}).
\subsection{Hypersurfaces and induced Spin$^c$ structures}
{\bf Spin$^c$ structures on manifolds:} Let $(M^n, g)$ be a Riemannian manifold of dimension $n\geqslant2$ without
boundary. We denote by $P_{\rm SO_n}M$ the 
${\rm SO}_n$-principal bundle over $M$ of positively oriented orthonormal frames. A
$\spinc$ structure of $M$ is a $\spin_n^c$-principal bundle $(P_{\spin_n^c} M,\pi,M)$
 and an $\Ss^1$-principal bundle $(P_{\Ss^1} M ,\pi,M)$ together with a double
covering given by  $\theta: P_{\spin_n^c} M \longrightarrow P_{\rm SO_n}M\times_{M}P_{\Ss^1} M$ such
that $\theta (ua) = \theta (u)\xi(a),$
for every $u \in P_{\spin_n^c} M$ and $a \in \spin_n^c$, where $\xi$ is the $2$-fold
covering of $\spin_n ^c$ over ${\rm SO}_n\times \Ss^1$. 
Let $\Sigma M := P_{\spin_n^c} M \times_{\rho_n} \Sigma_n$ be the associated spinor bundle
where $\Sigma_n = \C^{2^{[\frac n2]}}$ and $\rho_n : \spin_n^c
\longrightarrow  \eend(\Sigma_{n})$ denotes the complex spinor representation. A section of
$\Sigma M$ will be called a spinor field. The spinor bundle $\Sigma M$ is equipped with a
natural Hermitian scalar product denoted by $\<.,.\>$.\\\\
Additionally, any connection 1-form $A: T(P_{\Ss^1} M)\longrightarrow i\R$ on
$P_{\Ss^1} M$ and the connection 1-form 
$\omega^M$ on $ P_{\rm SO_n}M$, induce a connection
on the ${\rm SO}_n\times \Ss^1$-principal bundle $ P_{\rm SO_n}M\times_{M} P_{\Ss^1} M$, and hence 
a covariant derivative $\nabla$ on $\Gamma(\Sigma M)$ \cite{fr1,r2}. The curvature
of $A$ is an imaginary valued 2-form denoted by $F_A= dA$, i.e., $F_A = i\Omega$,
where $\Omega$ is a real valued 2-form on $P_{\Ss^1} M$. We know
 that $\Omega$ can be viewed as a real valued 2-form on $M$ $\cite{fr1, kn}$. In
this case $i\Omega$ is the curvature form of the auxiliary line bundle $L$. It is
the complex line bundle associated with the $\Ss^1$-principal bundle via the
standard representation of the unit circle.
For every spinor field $\psi$, the Dirac operator is locally defined 
by   $$D\psi =\sum_{j=1}^n e_j \cdot \nabla_{e_j} \psi,$$
where $\{e_1,\ldots,e_n\}$ is a local oriented orthonormal tangent frame and ``$\cdot$''  denotes the Clifford multiplication. The Dirac
operator is an elliptic, self-adjoint operator with respect to the $L^2$-scalar
product $(., .) = \int_M \<., .\> v_g $ and verifies, for any spinor field $\psi$, the Schr\"odinger-Lichnerowicz formula 
\begin{eqnarray}
D^2\psi=\nabla^*\nabla\psi+\frac{1}{4}S\psi+\frac{i}{2}\Omega\cdot\psi,
\label{bochner}
\end{eqnarray}
where $S$ is the scalar curvature of $M$,  $\nabla^*$ is the adjoint of $\nabla$ with respect to $(., .)$ and $\Omega\cdot$ is the extension of the Clifford multiplication to differential forms. For any $X \in \Gamma(TM)$, the Ricci identity is given by
\begin{eqnarray}
\sum_{k=1}^n e_k \cdot \mathcal{R}(e_k,X)\psi =
\frac 12 \Ric(X) \cdot \psi -\frac i2 (X\lrcorner\Omega)\cdot\psi,
\label{RRicci}
\end{eqnarray}
where $\Ric$ is the Ricci curvature of $(M^n, g)$ and $\mathcal{R}$ is the curvature tensor of the spinorial connection $\nabla$. In odd dimension, the volume form $\omega_{\C} := i^{[\frac{n+1}{2}]} e_1 \cdot...\cdot e_n$ acts on $\Sigma M$ as the identity, i.e., $\omega_\C \cdot\psi = \psi$ for any spinor $\psi \in \Gamma(\Sigma M)$. Besides, in even dimension, we have $\omega_\C^2 =1$. We denote by $\Sigma^\pm M$ the eigenbundles corresponding to the eigenvalues $\pm 1$, hence $\Sigma M = \Sigma^+ M \oplus \Sigma^- M$ and a  spinor field $\psi$ can be written $\psi = \psi^+ + \psi^-$. The conjugate $\overline \psi$ of $\psi$ is defined  by $\overline \psi = \psi^+ - \psi^-$.\\\\
Every spin manifold has a trivial $\spinc$ structure \cite{fr1}. In fact, we
choose the trivial line bundle with the trivial connection whose curvature $i\Omega$
is zero. Also every K\"ahler manifold $M$ of complex dimension $m$ ($n=2m$) has a canonical $\spinc$ structure coming from the complex structure $J$. Let $\ltimes$ be the K\"{a}hler form defined by the complex structure $J$, i.e. $\ltimes (X, Y)= g(JX, Y)$ for all vector fields $X,Y\in \Gamma(TM).$ The complexified tangent bundle $T^\C M =TM \otimes_\R \C$ decomposes into
$$T^\C M = T_{1,0} M\oplus T_{0,1} M,$$
where  $T_{1,0} M$ (resp. $T_{0,1} M$)  is the $i$-eigenbundle (resp. $-i$-eigenbundle) of the complex linear extension of the complex structure. Indeed,
$$T_{1,0}M = \overline{T_{0,1}M} = \{ X- iJX\ \ | X\in \Gamma(TM)\}.$$
Thus, the spinor bundle of the canonical $\spinc$ structure is given by $$\Sigma M = \Lambda^{0,*} M =\oplus_{r=0}^m \Lambda^r (T_{0,1}^* M),$$
where $T_{0,1}^* M$ is the dual space of  $T_{0,1} M$. The auxiliary bundle of this canonical $\spinc$ structure 
is given by  $L = (K_M)^{-1}= \Lambda^m (T_{0,1}^* M)$, where $K_M= \Lambda^m (T_{1,0}^* M)$ is the canonical bundle of $M$ \cite{fr1}. This line bundle $L$ has a canonical holomorphic connection
 induced from the Levi-Civita connection whose curvature form is given by $i\Omega = -i\rho$, where $\rho$ is the Ricci form given by $\rho(X, Y) = \Ric(JX, Y)$. Hence, this $\spinc$ structure carries  parallel spinors (the constant complex functions) lying in the set of complex functions $\Lambda^{0, 0}M \subset \Lambda^{0, *} M$  \cite{Moro1}.
Of course, we can define another $\spinc$ structure for which the spinor bundle is given by 
$\Lambda^{*, 0} M =\oplus_{r=0}^m \Lambda^r (T_{1, 0}^* M)$ and the auxiliary line bundle by $K_M$.  This $\spinc$ structure will be called the anti-canonical $\spinc$ structure \cite{fr1} and 
it carries  also parallel spinors (the constant complex functions) lying in the set of complex functions $\Lambda^{0, 0}M \subset \Lambda^{0, *} M$  \cite{Moro1}.\\\\
For any other $\spinc$ structure the spinorial bundle can be written as \cite{fr1, omu}: $$\Sigma M = \Lambda^{0,*}M\otimes\mathcal L,$$ where $\mathcal L^2 = K_M\otimes L$ and $L$  is the auxiliary bundle associated with this $\spinc$
structure. In this case, the $2$-form $\ltimes$ can be considered as an endomorphism of $\Sigma M$ via
 Clifford multiplication and 
it acts on a spinor field  $\psi$ locally by \cite{kirch, fr1}:
$$\ltimes\cdot\psi =  \frac 12 \sum_{j=1}^{m} e_j\cdot Je_j\cdot\psi.$$
Hence,  we have the well-known orthogonal splitting 
$$\Sigma M = \oplus_{r=0}^{m}\Sigma_rM,$$
where $\Sigma_rM$ denotes the eigensubbundle corresponding 
to the eigenvalue $i(m-2r)$ of $\ltimes$, with complex rank $\Big(^m_k\Big)$. The bundle $\Sigma_r M$ correspond 
to $\Lambda^{0, r}M\otimes\mathcal L$. Moreover,
$$\Sigma^+M = \bigoplus_{r \ \text{even}} \Sigma_r M\ \ \ \text{and} \ \ \ \Sigma^-M = \bigoplus_{r \ \text{odd}} \Sigma_rM.$$
For the canonical (resp. the anti-canonical) $\spinc$ structure, the subbundle $\Sigma_0M$ (resp. $\Sigma_mM$) is trivial, i.e., $\Sigma_0M = \Lambda^{0, 0}M \subset \Sigma^+M$ (resp. $\Sigma_mM = \Lambda^{0, 0}M$ which is in $\Sigma^+M$ if $m$ is even and in $\Sigma^-M$ if $m$ is odd). \\\\
{\bf Spin$^c$ hypersurfaces and the Gauss formula:} Let $N$ be an oriented ($n+1$)-dimensional Riemannian $\spinc$ manifold and $M
\subset N$ be an oriented hypersurface. The manifold $M$ inherits a $\spinc$
structure induced from the one on $N$, and we have \cite{r2}
$$ \Sigma M\simeq \left\{
\begin{array}{l}
\Sigma N_{|_M} \ \ \ \ \ \ \text{\ \ \ if\ $n$ is even,} \\\\
 \Sigma^+ N_{|_M}   \ \text{\ \ \ \ \ \ if\ $n$ is odd.}
\end{array}
\right.
$$
Moreover Clifford multiplication by a vector field $X$, tangent to $M$, is given by 
\begin{eqnarray}
X\bullet\phi = (X\cdot\nu\cdot \psi)_{|_M},
\label{Clifford}
\end{eqnarray}
where $\psi \in  \Gamma(\Sigma N)$ (or $\psi \in \Gamma(\Sigma^+ N)$ if $n$ is odd),
$\phi$ is the restriction of $\psi$ to $M$, ``$\cdot$'' is the Clifford
multiplication on $N$, ``$\bullet$'' that on $M$ and $\nu$ is the unit inner normal
vector. The connection 1-form defined on the restricted $\Ss^1$-principal bundle $(\PSM :=
\PSN_{|_M},\pi,M)$, is given by $A= {A^N}_{|_M} : T(\PSM) = T(\PSN)_{|_M}
\longrightarrow i\R.$ Then the curvature 2-form $i\Omega$ on the
$\Ss^1$-principal bundle $\PSM$ is given by $i\Omega= {i\Omega^N}_{|_M}$,
which can be viewed as an imaginary 2-form on $M$ and hence as the curvature form of
the line bundle $L$, the restriction of the auxiliary bundle $L^N$ to $M$. For every
$\psi \in \Gamma(\Sigma N)$ ($\psi \in \Gamma(\Sigma^+ N)$ if $n$ is odd), the real 2-forms
$\Omega$ and $\Omega^N$ are related by \cite{r2}
\begin{eqnarray}
(\Omega^N \cdot\psi)_{|_M} = \Omega\bullet\phi -
(\nu\lrcorner\Omega^N)\bullet\phi.
\label{glucose}
\end{eqnarray}
We denote by $\nabla^{\Sigma N}$ the spinorial Levi-Civita connection on $\Sigma
N$ and by $\nabla$ that on $\Sigma M$. For all $X\in \Gamma(TM)$, we have the spinorial Gauss formula \cite{r2}:
\begin{equation}
(\nabla^{\Sigma N}_X\psi)_{|_M} =  \nabla_X \phi + \frac 12 II X \bullet\phi,
\label{spingauss}
\end{equation}
where $II$ denotes the Weingarten map of the hypersurface. Moreover, Let $D^N$ and $D$ be the Dirac operators on $N$ and $M$, after
denoting by the same symbol any spinor and its restriction to $M$, we have
\begin{equation}
\widetilde D \phi = \frac{n}{2}H\phi -\nu\cdot D^N\phi-\nabla^{\Sigma N}_{\nu}\phi,
\label{diracgauss}
\end{equation}
where $H = \frac 1n \trace(II)$ denotes the mean curvature and $\widetilde D = D$ if $n$ is even and $\widetilde D =D \oplus(-D)$ if $n$ is odd.
\subsection{Basic facts about $\Ekt$ and their surfaces}
 We denote a 3-dimensional homogeneous manifolds with 4-dimensional isometry group by $\Ekt$. It is a Riemannian fibration over a simply connected 2-dimensional manifold $\M^2(\kappa)$ with constant curvature $\kappa$ and such that the fibers are geodesic. We denote by $\tau$ the bundle curvature, which measures the default of the fibration to be a Riemannian product. Precisely, we denote by $\xi$ a unit vertical vector field, that is tangent to the fibers. The vector field $\xi$ is a Killing field and satisfies for all  vector field $X$, 
$$\nablab_X\xi=\tau X\wedge\xi,$$ 
where $\nablab$ is the Levi-Civita connection and $\wedge$ is the exterior product. When $\tau$ vanishes, we get a product manifold $\M^2(\kappa)\times\R$. If $\tau\neq0$, these manifolds are of three types: They have the isometry group of the Berger spheres if $\kappa>0$, of the Heisenberg group $\mathrm{Nil}_3$ if $\kappa=0$ or of $\widetilde{\mathrm{PSL}_2(\R)}$ if $\kappa<0$.\\\\
Note that if $\tau =0$, then $\xi =\frac{\partial}{\partial t}$ is the unit vector field giving the orientation of $\R$ in the product $\M^2 (\kappa) \times \R$. The manifold $\Ekt$, with $\tau \neq 0$,  admits a local direct orthonormal frame $\{e_1, e_2, e_3\}$ with 
$$e_3 = \xi,$$
and such that the Christoffel symbols $\ov \Gamma^k_{ij} = \<\ov \nabla_{e_i}e_j, e_k\>$ are given by
\begin{eqnarray}\label{christoffel}
\left\lbrace  
\begin{array}{l}
\overline{\Gamma}_{12}^3=\overline{\Gamma}_{23}^1=-\overline{\Gamma}_{21}^3=-\overline{\Gamma}_{13}^2=\tau,\\ \\
\overline{\Gamma}_{32}^1=-\overline{\Gamma}_{31}^2=\tau-\frac{\kappa}{2\tau}, \\ \\
\overline{\Gamma}_{ii}^i=\overline{\Gamma}_{ij}^i=\overline{\Gamma}_{ji}^i=\overline{\Gamma}_{ii}^j=0,\quad\forall\,i,j\in\{1,2,3\},
\end{array}
\right. 
\end{eqnarray}
We call $\{e_1, e_2, e_3=\xi\}$ the canonical frame of $\Ekt$.\\\\
Let $M$ be a simply connected orientable  surface of $\Ekt$ with shape operator $II$ associated with the unit inner normal vector $\nu$. Moreover, we denote $\xi=T+f\nu$ where the function $f$ is the normal component of $\xi$ and $T$ is its tangential part. We introduce the following notion of compatibility equations.
\begin{defn}[{\bf Compatibility equations}]\label{comp}
We say that $(M,\langle.,.\rangle,E,T,f)$ satisfies the compatibility equations for $\Ekt$ if and only if for any $X,Y,Z\in \Gamma(TM)$,
\begin{align}\label{gauss}
K=\ddet(E)+\tau^2+(\kappa-4\tau^2)f^2
\end{align}
\beqt\label{codazzi}
\nabla_XEY-\nabla_YEX-W[X,Y]=(\kappa-4\tau^2) f(\langle Y,T\rangle X-\langle X,T\rangle Y),
\eeqt
\beqt\label{cond1}
\nabla_XT=f(EX-\tau JX),
\eeqt
\beqt\label{cond2}
df(X)=-\langle EX-\tau JX,T\rangle,
\eeqt
where $K$ is the Gauss curvature of $M$.
\end{defn}
\begin{rem}
The relations (\ref{gauss}) and (\ref{codazzi}) are the Gauss and Codazzi equations for an isometric immersion into $\Ekt$ obtained by a computation of the curvature tensor of  $\Ekt$. Equations \eqref{cond1} and \eqref{cond2} are coming from the fact that $\nablab_X\xi=\tau X\wedge\xi$.
\end{rem}
In \cite{Da,Da2}, Daniel proves that these compatibility equations are necessary and sufficient for the existence of an isometric immersion $F$ from $M$ into $\Ekt$  with shape operator $dF\circ E\circ dF^{-1}$ and so that $\xi=dF(T)+f\nu$.

\subsection{Basic facts about $\CM$ and their real hypersurfaces}
\label{cm}
Let $(\CM, J, \overline{g})$ be the complex space form of constant holomorphic sectional curvature $4c \neq 0$ and complex dimension $2$, that is for $c=1$, $\CM$ is the complex projective space $\CP$ and if $c=-1$, $\CM$ is the complex hyperbolic space $\CH$. It is a well-known fact that the curvature tensor $\overline{R}$ of $\CM$ is given by
\beQ \overline{g}\big(\overline{R}(X,Y)Z,W\big)&=&c\Big\{ \overline{g}(Y,Z)\overline{g}(X,W)-\overline{g}(X,Z)\overline{g}(Y,W)+\overline{g}(JY,Z)\overline{g}(JX,W)\\
&&-\overline{g}(JX,Z)\overline{g}(JY,W)-2\overline{g}(JX,Y)\overline{g}(JZ,W)\Big\},
\eeQ
for all $X,Y,Z$ and $W$ tangent vector fields to $\CM$.\\\\
Let $M^3$ be an oriented real hypersurface of $\CM$ endowed with the metric $g$ induced by $\overline{g}$. 
We denote by $\nu$ a normal unit inner  vector globally defined on $M$ and  by $II$ the shape operator of this 
immersion. Moreover, the complex structure $J$ induces on $M$ an almost contact metric structure $(\Chi, \xi, \eta, g)$, 
where $\Chi$ is the $(1,1)$-tensor defined by $g(\Chi X,Y)=\overline{g}(JX,Y)$ for all $X,Y\in \Gamma(TM)$, $\xi=-J\nu$ 
is a tangent vector field and $\eta$ the $1$-form associated with $\xi$, that is so that 
$\eta(X)=g(\xi,X)$ for all $X\in\Gamma(TM)$. Then, we see easily that the following holds:
\beqt
\Chi^2X=-X+\eta(X)\xi,\quad g(\xi,\xi)=1,\quad\text{and}\quad \Chi\xi=0.
\eeqt
Here, we recall that given an almost contact metric structure $(\Chi, \xi, \eta, g)$ one defines a $2$-form 
$\varTheta$ by $\varTheta(X, Y) = g(X, \Chi Y)$ for all $X, Y\in \Gamma(TM)$. Now, $(\Chi, \xi, \eta, g)$ is said to satisfy 
the contact condition if $-2 \varTheta =d\eta$ and if it is the case, $(\Chi, \xi, \eta, g)$ is called a contact metric 
structure on $M$. A contact metric structure $(\Chi, \xi, \eta, g)$  is called a Sasakian structure (and $M$ a Sasaki manifold) if $\xi$ is a Killing vector field (or equivalently, $\Chi = \nabla \xi$) and 
$$(\nabla_X\Chi )Y = \eta (Y) X - g(X, Y) \xi, \ \ \text{for all} \ \ X, Y \in \Gamma(TM).$$
From the relation between the Riemannian connections of $\CM$ and $M$, $\overline{\nabla}_XY=\nabla_XY+g(IIX,Y)\nu$, we deduce the two following identities:
\beqt\label{condd3}
(\nabla_X\Chi)Y=\eta(Y)IIX-g(IIX,Y)\xi,
\eeqt
\beqt\label{condd4}
\nabla_X\xi=\Chi IIX.
\eeqt
From the expression of the curvature of $\CM$ given above, we deduce the Gauss and Codazzi equations. First, the Gauss equation says that for all $X,Y,Z,W\in \Gamma(TM)$,
\beq\label{gaussCMM}\nonumber
g(R(X,Y)Z,W)&=&c\Big\{ g(Y,Z)\overline{g}(X,W)-g(X,Z)g(Y,W)+g(\Chi Y,Z)g(\Chi X,W)\\ 
&&-g(\Chi X,Z)g(\Chi Y,W)-2g(\Chi X,Y)g(\Chi Z,W)\Big\}\\
&&+g(IIY,Z)g(IIX,W)-g(IIX,Z)g(IIY,W).\nonumber
\eeq
The Codazzi equation is 
\beq\label{codazziCMM}
d^{\nabla}II(X,Y)=c\big(\eta(X)\Chi Y-\eta(Y)\Chi X-2g(\Chi X,Y)\xi\big).
\eeq
Now, we ask if the Gauss equation \eqref{gaussCMM} and the Codazzi equation \eqref{codazziCMM} are sufficient to get an isometric immersion of $(M,g)$ into $\CM$.
\begin{defn}[{\bf Compatibility equations}]\label{compp}
Let $(M^3,g)$ be a simply connected oriented Riemannian manifold endowed with an almost contact metric structure $(\Chi, \xi, \eta)$ and $E$ be a field of symmetric endomorphisms on $M$. We say that $(M, g, E, \Chi, \xi, \eta)$ satisfies the compatibility equations for $\CM$ if and only if for any $X,Y,Z,W\in\Gamma(TM)$, we have
\begin{eqnarray}
\label{gaussCM}\nonumber
g(R(X,Y)Z,W)&=&c\Big\{ g(Y,Z)\overline{g}(X,W)-g(X,Z)g(Y,W)+g(\Chi Y,Z)g(\Chi X,W)\\ 
&&-g(\Chi X,Z)g(\Chi Y,W)-2g(\Chi X,Y)g(\Chi Z,W)\Big\}\\
&&+g(EY,Z)g(EX,W)-g(EX,Z)g(EY,W),\nonumber
\end{eqnarray}
\begin{eqnarray}\label{codazziCM}
d^{\nabla}E(X,Y)=c\big(\eta(X)\Chi Y-\eta(Y)\Chi X-2g(\Chi X,Y)\xi\big).
\end{eqnarray}
\begin{eqnarray}
\label{cond3}
(\nabla_X\Chi)Y=\eta(Y)EX-g(EX,Y)\xi,
\end{eqnarray}
\begin{eqnarray}
\label{cond4}
\nabla_X\xi=\Chi EX.
\end{eqnarray}
\end{defn}
In \cite{PT}, P. Piccione and D. V. Tausk proves that the Gauss equation  \eqref{gaussCM} and the Codazzi equation \eqref{codazziCM} together with (\ref{cond3}) and (\ref{cond4}) are necessary and sufficient for the existence of an isometric immersion from $M$ into $\CM$ such that the complex structure of $\CM$ over $M$ is given by $J= \Chi  + \eta(\cdot) \nu$.

\section{Isometric immersions into $\Ekt$ via spinors}
The manifold $\Ekt$ has a $\spinc$ structure carrying a Killing spinor with Killing constant $\frac{\tau}{2}$. The restriction of this $\spinc$ structure to any surface $M$ defines a $\spinc$ structure on $M$ with a special  spinor field. This spinor field characterizes the isometric immersion of $M$ into $\Ekt$. 
\subsection{Special spinors fields on $\Ekt$ and their surfaces}
On $\spinc$ manifolds, A. Moroianu defined projectable spinors for arbitrary Riemannian submersions of $\spinc$ manifolds with $1$-dimensional totally geodesic fibers \cite{Morothese, Moro2}. These spinors will be used to get a Killing spinor on $\Ekt$.
\begin{prop}
The canonical $\spinc$ structure on $\M^2 (\kappa)$ induces a $\spinc$ structure on $\Ekt$ carrying 
a Killing spinor with Killing constant $\frac{\tau}{2}$.
\end{prop}
{\bf Proof}: By enlargement of the group structures, the two-fold covering $\theta: P_{\spin_2^c} \M \longrightarrow P_{{\rm SO}_2} \M\times_{\M}P_{\Ss^1} \M$, gives a two-fold covering 
$$\theta: P_{\spin_3^c} \M \longrightarrow P_{{\rm SO}_3} \M\times_{\M}P_{\Ss^1} \M,$$
which, by pull-back through $\pi$, gives rise to a $\spinc$ structure on $\overline M:= \Ekt$ \cite{Moro2, Morothese} and the following diagram commutes
$$
\xymatrix{
P_{\spin_3^c} \overline M \ar[d]^{\pi^*\theta} \ar[r] & P_{\spin_3^c} \M\ar[d]_\theta 
& & \\
P_{\rm SO_3} \overline M \times_{\overline M} P_{\Ss^1} \overline M   \ar[r] &  P_{{\rm SO}_3} \M\times_{\M}P_{\Ss^1} \M } 
$$
The next step is to relate the covariant derivatives of spinors on $\M$ and $\overline M$. We point out an important detail: Since we are actually interested to get a Killing spinor on $\overline M$, the connection on $P_{\Ss^1}\overline M$ (which defines the covariant derivative of spinors on $\overline M$) that we will consider will be the pull-back connection if $\tau =0$ and  will not be the pull-back connection if $\tau \neq 0$. Hence, when $\tau =0$, the connection $A_0$ on $P_{\Ss^1} \overline M$ is given by
$$A_0((\pi^* s)_* (X^*)) = A(s_* X)\ \ \ \ \ \text{and} \ \ \ \ \ A_0((\pi^* s)_* \xi) = 0.$$
Now, if $\tau \neq 0$, we consider a connection $A_0$ on $P_{\Ss^1} \overline M$  given by
$$A_0((\pi^* s)_* (X^*)) = A(s_* X)\ \ \ \ \ \text{and} \ \ \ \ \ A_0((\pi^* s)_* \xi) = -i(2\tau-\frac{\kappa}{2\tau}),$$
where $e_3 = \xi$ is the vertical vector field on $\Ekt$ if $\tau \neq 0$ or $e_3 = \partial t$ if $\tau =0$, $X^*$ is the horizontal left of a vector field $X$ on $\M$, $A$ is the connection defined on $P_{\Ss^1 }\M$ and $s$ a local section of $P_{\Ss^1 }\M$. Recall that we have an identification of the pull back $\pi^* \Sigma \M$ with $\Sigma \overline  M$ \cite{Moro2, Morothese}, and with respect to this identification, if $X$ is a vector field and $\psi$ a spinor field on $\M$, then
\begin{eqnarray}
 X^* \cdot \pi^*\psi = \pi^* (X\cdot\psi)\ \ \ \ \ \text{and}\ \ \ \ \ \ \xi\cdot\pi^* \psi = -i\pi^*(\overline \psi).
\label{pro}
\end{eqnarray}
The sections of $\Sigma\overline M$ which can be written as pull-back of sections of $\Sigma \M$ are called projectable spinors \cite{Moro2, Morothese}. Now, we relate the covariant derivative $\nabla^{\Ekt}$ of projectable spinors on $\Ekt$  to the covariant derivative $\nabla$ of spinors on $\M$.
In fact, any spinor field $\psi$ is locally written as  $\psi = [\widetilde{b \times s},\sigma]$,
where 
$b = (e_1, e_2)$ is a base of $\M^2(\kappa)$, $s: U\longrightarrow 
P_{\Ss^1} \M$ is a local section of $P_{\Ss^1}\M$ and $\widetilde{b \times s}$ is the lift of the local section $b \times s: U \rightarrow P_{SO_2}\M\times_{\M} P_{\Ss^1}\M$ by the 2-fold covering. Then $\pi^*\psi$ can be expressed as $\pi^*\psi = [\pi^* (\widetilde{b \times s}), \pi^*\sigma]$. It is easy to see that the projection $\pi^*(\widetilde{b \times s})$ onto $P_{SO_3} \overline M$ is the canonical frame $(e_1^*, e_2^* , e_3=\xi)$ and its projection onto $P_{\Ss^1}\overline M$ is just $\pi^* \sigma$. We have
\begin{eqnarray*}
 \nabla^{\Ekt}_{e_1^*} \pi^*\psi &=& [\pi^* (\widetilde{b \times s}), e_1^* (\pi^* \sigma)] + \frac 12  \overline g(\overline \nabla_{e_1^*} e_1^*, e^*_2)e^*_1\cdot e^*_2\cdot\pi^*\psi \\&& + \frac 12 \sum_{j=1}^2 \overline g(\overline\nabla_{e_1^*}e_j^*, e_3) e_j^*\cdot e_3 \cdot \pi^*\psi + \frac 12 A_0((\pi^*s)_* e_1^*) \pi^*\psi\\
 &\stackrel{(\ref{christoffel})}{=}&    [\pi^* (\widetilde{b \times s}), \pi^* (e_1(\sigma))] + \frac 12 g(\nabla_{e_1}e_1, e_2) \pi^* (e_1\cdot e_2\cdot\psi)\\ && + \frac {\tau}{2} e_2^*\cdot e_3\cdot\pi^*\psi  + \frac12 A(s_*X) \pi^*\psi\\&=&
\pi^* \Big([(\widetilde{b \times s}), (e_1(\sigma))]   + \frac 12  g(\nabla_{e_1}e_1, e_2) e_1 \cdot e_2 \cdot\psi \\ &&+ \frac{\tau}{2} e_1^*\cdot\psi +\frac12 A(s_*X) \psi\Big) \\&=&
 \pi^*(\nabla_{e_1}\psi) + \frac{\tau}{2} e_1\cdot\pi^*\psi.\\
\end{eqnarray*}
The same holds for $e_2^*$. Similary, if $\tau  \neq 0$ we have
\begin{eqnarray*}
 \nabla^{\Ekt}_{e_3} \pi^*\psi &=& [\pi^* (\widetilde{b \times s}), e_3 (\pi^* \sigma)] + \frac 12  \overline g(\overline \nabla_{e_3} e_1^*, e^*_2)e^*_1\cdot e^*_2\cdot\pi^*\psi \\&& + \frac 12 \sum_{j=1}^2 \overline g(\nabla_{e_3}e_j^*, e_3) e_j^*\cdot e_3 \cdot \pi^*\psi + \frac 12 A_0((\pi^*s)_* e_3) \pi^*\psi\\
 &\stackrel{(\ref{christoffel})}{=}& \frac{1}{2}\Big({\frac{\kappa}{2\tau}-\tau}\Big) e_1^*\cdot e_2^*\cdot\pi^*\psi - \frac{i}{2}\Big( 2\tau -\frac{\kappa}{2\tau}\Big) \pi^*\psi\\
&=& \frac{1}{2}\Big({\frac{\kappa}{2\tau}-\tau}\Big) e_3\cdot\pi^*\psi + \frac{1}{2}\Big( 2\tau -\frac{\kappa}{2\tau}\Big) e_3\cdot\pi^*\overline\psi.
\end{eqnarray*}
Now, the canonical $\spinc$ structure on $\M^2(\kappa)$ carries a parallel spinor  $\psi \in \Gamma(\Sigma_0 \M) \subset \Gamma(\Sigma^+ \M)$, so $\overline \psi = \psi.$ For this canonical $\spinc$ structure, the determinant line bundle corresponding to $P_{\Ss^1} \M$ is $K_{\M}^{-1}$ and the connection $1$-form $A$ on  $P_{\Ss^1} \M$ is the connection for the Levi-Civita connection extended to $K_\M^{-1}$.  Hence, the spinor $\pi^*\psi$ is a Killing spinor field on $\Ekt$, because
$$ \nabla^{\Ekt}_{e_j^*} \pi^*\psi = \frac{\tau}{2} e_j^*\cdot\pi^*(\psi), \ \ \ \text{for}\ \ \  j=1, 2\ \ \text{and}\ \ \  \nabla^{\Ekt}_{\xi} \pi^*\psi = \frac{\tau}{2} \xi\cdot\pi^*\psi.$$
Now, if $\tau=0$, a same computation of $\nabla^{\Ekt}_{e_3} \pi^*\psi$ gives that $\pi^*\psi$  is a parallel spinor field 
 on $\Ekt$.
\begin{rem}
Every Sasakian manifold has a canonical $\spinc$ structure: In fact, giving a Sasakian structure on a manifold $(M^n, g)$ is equivalent to give a K\"{a}hler structure on the cone  over $M$. The cone over $M$ is the manifold $ M\times_{r^2} \R^+$ equipped with the metric $r^2 g + dr^2$. Moreover, there is a $1$-$1$-correspondence between $\spinc$ structures on $M$ and that on its cone \cite{Moro1}.  Hence, every Sasakian manifold has a canonical (resp. anti-canonical) $\spinc$ structure coming from the canonical one (resp. anti-canonical one) on its cone. \\\\
In \cite{Moro1}, A. Moroianu classified all complete simply connected $\spinc$ manifolds carrying real  Killing spinors and  he proved that the only complete simply connected $\spinc$ manifolds carrying real Killing spinors (other than the $Spin$ manifolds) are the non-Einstein Sasakian manifolds endowed with their canonical (or anti-canonical) $\spinc$ structure. \\\\
The manifold $\Ekt$ is a complete simply connected non-Einstein manifold and hence the only $\spinc$ structure carrying a Killing spinor is the canonical one (or the anti-canonical). Hence, the $\spinc$ structure on $\Ekt$ described above, (i.e. the one coming from $\M^2(\kappa)$) is nothing than the canonical $\spinc$ structure coming from the Sasakian structure.\\\\
We point out that, in a similar way, the anti-canonical $\spinc$ structure on $\M^2(\kappa)$
 (carrying a parallel spinor field lying in $\Sigma^-\M$) induces also on $\Ekt$ the anti-canonical $\spinc$ structure with 
a Killing spinor $\pi^*\psi$ of Killing constant $\frac{\tau}{2}$ if $\tau \neq 0$ and a parallel spinor $\pi^*\psi$ if 
$\tau =0$. In both cases, we have $\xi\cdot \pi^*\psi = -i \pi^*\overline \psi =  i\pi^*\psi$. For $\tau \neq 0$, the connection $A_0$ is chosen to be 
$$A_0((\pi^* s)_* (X^*)) = A(s_* X)\ \ \ \ \ \text{and} \ \ \ \ \ A_0((\pi^* s)_* \xi) = i(2\tau-\frac{\kappa}{2\tau}).$$
When $\tau =0$, it is the pull-back connection.
\end{rem}
From now, we will denote the Killing spinor field $\pi^*\psi$ on $\Ekt$ by $\psi$. Since, it is a Killing spinor, we have 
$$(\nabla^{\Ekt})^*\nabla^{\Ekt} \psi = \frac{3\tau^2}{4} \psi\ \ \ \ \ \text{and}\ \ \ \ \ \ D^{\Ekt}\psi = -\frac{3\tau}{2} \psi.$$
By the Schr\"{o}dinger-Lichnerowicz formula, we get $$\frac i2 \Omega^{\Ekt} \cdot\psi = \frac{3 \tau^2}{2} \psi -\frac{(\kappa-\tau^2)}{2}\psi,$$
where $i\Omega^{\Ekt}$ is the curvature 2-form of the auxiliary line bundle associated with the $\spinc$ structure. Finally, 
\begin{eqnarray}
 \Omega^{\Ekt}\cdot\psi = i(\kappa-4\tau^2) \psi.
\label{omega}
\end{eqnarray}
\subsection{Spinorial characterization of surfaces of $\Ekt$}
Let $\kappa,\tau\in\R$ with $\kappa-4\tau^2\neq0$ and $M$ be a simply connected oriented Riemannian surface immersed into $\Ekt$. The vertical vector field $\xi$ is written  $\xi = T +f\nu$ where $T$ be a vector field on $M$ and $f$ a real-valued function on $M$ so that $f^2+||T||^2=1$. We endowed $\Ekt$ with the $\spinc$ structure described above, carrying a Killing spinor of  Killing constant $\frac{\tau}{2}$.
\begin{lem}
 The restriction $\varphi$ of the Killing spinor $\psi$ on $\Ekt$ is a solution of the following equation 
\begin{eqnarray}
 \nabla_X\varphi + \frac 12 II X\bullet\varphi - i\frac{\tau}{2} X\bullet\ov\varphi=0,
\label{restricted}
\end{eqnarray}
called the restricted Killing spinor equation. Moreover, $f= \frac{<\varphi, \overline \varphi>}{\vert\varphi\vert^2}$ and the curvature 2-form of the connection on the auxiliary line bundle associated with the induced $\spinc$ structure is given by $\Omega(t_1, t_2)= -(\kappa-4\tau^2)f$, in any local orthonormal frame $\{t_1, t_2\}$.
\label{lemEkt}
\end{lem}
{\bf Proof:} We restrict the $\spinc$ structure on $\Ekt$ to $M$. By the Gauss formula (\ref{spingauss}), the restriction $\varphi$ of the Killing spinor $\psi$ on $\Ekt$ satisfies
$$\nabla_X\varphi + \frac 12 II X \bullet\varphi - \frac{\tau}{2} {X\cdot\psi}_{\vert_M}=0.$$
Let $\{t_1, t_2, \nu\}$ be a local orthonormal frame of $\Ekt$ such that $\{t_1, t_2\}$ is a local orthonormal frame of $M$ and $\nu$  a unit normal vector field of the surface. The action of the volume forms on $M$ and $\Ekt$ gives
\begin{eqnarray*}
 X\bullet\overline\varphi &=&  i(X\bullet t_1 \bullet t_2 \bullet \varphi) \\&=& i(X\cdot\nu\cdot t_1 \cdot t_2 \cdot \psi)_{\vert_M} \\&=& -i (X\cdot\psi)_{\vert_M}, 
\end{eqnarray*}
which gives Equation (\ref{restricted}). The vector field $T$ splits into $T=\nu_1+h \xi$ where $\nu_1$ is a vector field generated by $e_1$ and $e_2$ and $h$ a real function. The scalar product of $T$  by $\xi = T+f\nu$ and the scalar product of $T =\nu_1 + h \xi$ by $\xi$  gives  $||T||^2=h$ which means that $h=1-f^2$. Hence, the normal vector field $\nu$ can be written as $\nu=f \xi-\frac{1}{f}\nu_1.$ As we mentionned before, the $\spinc$ structure on $\Ekt$ induces a $\spinc$ structure on $M$ with induced auxiliary line bundle. Next, we want to prove that the curvature 2-form of the connection on the auxiliary line bundle of $M$ is equal to $i\Omega(t_1,t_2)=-i (\kappa-4\tau^2)f$. Since the spinor $\psi$ is Killing, the equality (\ref{RRicci}) gives,  for all $X\in T(\Ekt)$ 
\begin{eqnarray}
\mathrm {Ric}^{\Ekt}(X)\cdot\psi - i (X\lrcorner \Omega^{\Ekt})\cdot\psi = 2\tau^2 X\cdot\psi,
\label{RRRicci}
\end{eqnarray}
Where $\Ric$ is the Ricci tensor of $\Ekt$. Therfore, we compute,
\begin{eqnarray*}
(\nu\lrcorner \Omega^{\Ekt})\bullet\varphi &=&(\nu\lrcorner \Omega^{\Ekt})\cdot\nu\cdot\psi|_M \\ &=&
i (2\tau^2 \psi + \nu \cdot\mathrm{Ric}^{\Ekt}\ \nu\cdot\psi)_{\vert_M}.
\end{eqnarray*}
But  we have $\Ric^{\Ekt}e_3 = 2\tau^2 e_3$, $\Ric^{\Ekt}e_1 = (\kappa-2\tau^2) e_1$ and $\Ric^{\Ekt}e_2 =(\kappa-2\tau^2) e_2$. Hence,
\begin{eqnarray*}
 \Ric^{\Ekt}\nu &=& f \Ric^{\Ekt} e_3 -\frac 1f \Ric^{\Ekt} \nu_1 = 2\tau^2 f e_3 -\frac 1f(\kappa-2\tau^2) \nu_1\\
&=& 2\tau^2 f e_3 + (\kappa-2\tau^2) (\nu -f e_3)\\
&=& -(\kappa-4\tau^2) f e_3 +(\kappa-2\tau^2)\nu.
\end{eqnarray*}
We conclude using Equation (\ref{pro}) that
$$(\nu\lrcorner \Omega^{\Ekt})\bullet\varphi = -i(\kappa-4\tau^2)\varphi - (\kappa-4\tau^2)f (\nu\cdot\psi)_{\vert_M}.$$
By Equation (\ref{glucose}), we get that $\Omega\bullet \varphi = -(\kappa-4\tau^2) f (\nu\cdot\psi)_{\vert_M}$. The scalar product of the last equality with $t_1 \bullet t_2 \bullet \varphi$ gives 
$$\Omega(t_1,t_2)|\varphi|^2=f (\kappa-4\tau^2) (\psi,t_1\cdot t_2\cdot\nu\cdot \psi)|_M= -f (k-4\tau^2)|\varphi|^2.$$
We write in the frame $\{t_1,t_2,\nu\}$ 
\begin{equation}\label{eq:3}
\Omega^{\Ekt}(t_1,t_2)t_1\cdot t_2\cdot\psi+\Omega^{\Ekt}(t_1,\nu)t_1\cdot\nu\cdot\psi+\Omega^{\Ekt}(t_2,\nu)t_2\cdot\nu\cdot\psi=i(\kappa-4\tau^2)\psi.
\end{equation}
But we know that $\Omega^{\Ekt}(t_1,t_2) =\Omega (t_1,t_2) =-(\kappa-4\tau^2)f$. For the other terms, we compute 
$$
\Omega^{\Ekt}(t_1,\nu)=\Omega^{\Ekt}(t_1,\frac{1}{f} e_3-\frac{1}{f} T)=- \frac{1}{f}g(T,t_2)\Omega^{\Ekt}(t_1,t_2)= (\kappa-4\tau^2)g(T,t_2),
$$
where the term $\Omega^{\Ekt}(t_1, e_3)$ vanishes since by Equation (\ref{RRRicci}) we have $e_3 \lrcorner \Omega^{\Ekt} =0$. Similarly, we find that $\Omega^{\Ekt}(t_2,\nu)=-(\kappa-4\tau^2)g(T,t_1).$ By substituting these values into \eqref{eq:3} and taking Clifford multiplication with $t_1\cdot t_2$, we get 
$$T\bullet\varphi =- f\varphi +\overline\varphi.$$
Finally, take the real part of the scalar product of the last equation by $\varphi$, we get $f= \frac{<\varphi, \overline \varphi>}{\vert\varphi\vert^2}$.
\begin{rem}
Using also the Equation $T\bullet\varphi =- f\varphi +\overline\varphi$, we can deduce that
$$g(T,t_1)=\Re\<it_2\bullet\varphi,\frac{\varphi}{|\varphi|^2}\> \quad\text{and}\quad g(T,t_2)=-\Re\<it_1\bullet\varphi,\frac{\varphi}{|\varphi|^2}\>.$$
\end{rem}
\begin{prop}
 Let $(M^2,g)$ be an oriented $\spinc$ surface carrying  a non-trivial solution $\varphi$ of the following equation
$$ \nabla_X\varphi + \frac 12 EX\bullet\varphi - i\frac{\tau}{2} X\bullet\ov\varphi=0,$$
where $E$ denotes a symmetric tensor field defined on $M$. Moreover, assume that the curvature 2-form of the asssociated auxiliary bundle satisfies $i\Omega(t_1, t_2)=  -(\kappa-4\tau^2) f =-(\kappa-4\tau^2)  \frac{<\varphi, \overline \varphi>}{\vert\varphi\vert^2}$ in any local orthonormal frame $\{t_1, t_2\}$ of $M$. Then, there exists an isometric immersion of $(M^2, g)$ into $\Ekt$ with shape operator $E$, mean curvature $H$ and such that, over $M$, the vertical vector is $\xi = dF(T) + f\nu$, where $\nu$ is the unit normal vector to the surface and $T$ is the tangential part of $\xi$ given by 
$$g(T,t_1)=\Re\<it_2\bullet\varphi,\frac{\varphi}{|\varphi|^2}\> \quad\text{and}\quad g(T,t_2)=-\Re\<it_1\bullet\varphi,\frac{\varphi}{|\varphi|^2}\>.$$

\label{propEkt}
\end{prop}

{\bf Proof:} We compute the action of the spinorial curvature tensor $\mathcal{R}$ on $\varphi$. We have
\begin{eqnarray*}
 \nabla_{t_1}\nabla_{t_2} \varphi & =& -\frac 12 \nabla_{t_1} E(t_2)\bullet\varphi +\frac 14 E(t_2)\bullet E(t_1)\bullet\varphi - \frac{\tau}{4} E(t_2)\bullet t_2\bullet\varphi \\ && -\frac{\tau}{2} \nabla_{t_1}(t_1)\bullet\varphi + \frac{\tau}{4}t_1\bullet E(t_1)\bullet\varphi-\frac{\tau^2}{4}t_1\bullet t_2 \bullet\varphi.
\end{eqnarray*}
As well as
\begin{eqnarray*}
 \nabla_{t_2}\nabla_{t_1} \varphi &=& -\frac 12\nabla_{t_2} E(t_1)\bullet\varphi + \frac 14 E(t_1)\bullet E(t_2)\bullet\varphi - \frac{\tau}{4} E(t_1)\bullet t_1\bullet\varphi \\&&-\frac{\tau}{2} \nabla_{t_2}{t_2}\bullet\varphi + \frac{\tau}{4}t_2\bullet E(t_2)\bullet\varphi-\frac{\tau^2}{4}t_2\bullet t_1 \bullet\varphi.
\end{eqnarray*}
So, taking into account that $[t_1, t_2] = \nabla_{t_1}t_2 -\nabla_{t_2}t_1$, a straightforward computation gives 
$$ \mathcal{R}(t_1, t_2)\varphi = -\frac 12 (d^\nabla E)(t_1,t_2)\bullet\varphi - \frac 12 {\rm det}\,E \ t_1\bullet t_2\bullet\varphi - \frac{\tau^2}{2}t_1 \bullet t_2 \bullet \varphi.$$
On the other hand, it is well known that 
$$\mathcal {R}(t_1, t_2)\varphi = - \frac 12 R_{1212}\ t_1\bullet t_2 \bullet \varphi +\frac{i}{2} \Omega(t_1, t_2) \varphi.$$
Therefore, we have
\begin{eqnarray}
 (R_{1212} -\mathrm{det} E - \tau^2)t_1\bullet t_2 \bullet\varphi = (d^\nabla E(t_1,t_2)-if(\kappa-4\tau^2))\varphi.
\label{GC}
\end{eqnarray}
Now, let $T$ a vector field of $M$ given by 
$$g(T,t_1)|\varphi|^2=\Re\<it_2\bullet\varphi,\varphi\> \quad\text{and}\quad g(T,t_2)|\varphi|^2=-\Re\<it_1\bullet\varphi,\varphi\>.$$
It is easy to check that $T\bullet\varphi = -f\varphi + \overline \varphi$ and hence $f^2 + \Vert T\Vert^2 =1$. In the following, we will prove that the spinor field $\theta:=i\varphi-if\overline{\varphi}+JT\bullet\varphi$ is zero. For this, it is sufficient to prove that its norm vanishes. Indeed, we compute 
\begin{equation}
|\theta|^2=|\varphi|^2+f^2|\varphi|^2+||T||^2|\varphi|^2-2\Re\<i\varphi,if\ov{\varphi}\>+2\Re\<i\varphi,JT\bullet\varphi\>
\label{eq:5}
\end{equation}
Therefore, Equation \eqref{eq:5} becomes 
\begin{eqnarray*}
|\theta|^2&=&2|\varphi|^2-2f^2|\varphi|^2+2\Re\<i\varphi,JT\bullet\varphi\>\\
&=&2|\varphi|^2-2f^2|\varphi|^2+2 g(JT,t_1)\Re\<i\varphi,t_1\bullet\varphi\>+2g(JT,t_2) \Re\<i\varphi,t_2\bullet\varphi\>\\
&=&2|\varphi|^2-2f^2|\varphi|^2+2 g(JT,t_1) g(T,t_2)|\varphi|^2-2g(JT,t_2) g(T,t_1)|\varphi|^2\\
&=&2|\varphi|^2-2f^2|\varphi|^2-2g(T,t_2)^2|\varphi|^2-2g(T,t_1)^2|\varphi|^2\\
&=&2|\varphi|^2-2f^2|\varphi|^2-2||T||^2|\varphi|^2=0.
\end{eqnarray*}
Thus, we deduce $if\varphi=-f^2t_1\bullet t_2\bullet\varphi-fJT\bullet\varphi$, where we use the fact that $\overline\varphi=i t_1\bullet t_2\bullet\varphi$. In this case, Equation (\ref{GC}) can be written as 
$$(R_{1212} -\mathrm{det} E - \tau^2 - (\kappa-4\tau^2)f^2)t_1\bullet t_2 \bullet\varphi = (d^\nabla E(t_1,t_2) + (\kappa-4\tau^2)JT)\bullet\varphi.$$
This is equivalent to say that both terms $R_{1212} -\mathrm{det} E - \tau^2 - (\kappa-4\tau^2)f^2$ and $d^\nabla E(t_1,t_2) + (\kappa-4\tau^2)JT$ are equal to zero.  In fact, these are the Gauss-Codazzi equations in Definition \ref{comp}. In order to obtain the two other equations, we simply compute the derivative of $T\bullet\varphi=-f\varphi+\overline\varphi$ in the direction of $X$ in two ways. First, using that $iX\bullet\overline\varphi = JX\bullet\varphi$, we have 
\begin{eqnarray}
\nabla_X T \bullet \varphi+T\bullet\nabla_X\varphi&=&\nabla_X T \bullet \varphi-\frac{1}{2}T\bullet EX\bullet\varphi +i\frac {\tau}{2} T\bullet X\bullet\overline\varphi\nonumber\\
&=&\nabla_X T \bullet \varphi-\frac{1}{2}T\bullet EX\bullet\varphi +\frac {\tau}{2} T\bullet JX\bullet\varphi.
\label{c1}
\end{eqnarray}
On the other hand, we have 
\begin{eqnarray}
\nabla_X (T\bullet\varphi)&=& -X(f)\varphi-f\nabla_X\varphi+\nabla_X\overline\varphi
\nonumber\\ &=&-X(f)\varphi+\frac{1}{2}f EX\bullet\varphi+\frac{1}{2} EX\bullet\overline{\varphi} -i\frac{\tau}{2} f X\bullet\overline\varphi -\frac i2 \tau X\bullet\varphi\nonumber\\
&=& -X(f)\varphi+\frac{1}{2}f EX\bullet\varphi -\frac 12 f\tau JX\bullet\varphi\nonumber\\ &&
+\frac 12 EX\bullet(T\bullet\varphi+f\varphi) - \frac{i}{2} \tau X\bullet\varphi\nonumber\\
&=&-X(f)\varphi+\frac{1}{2}f EX\bullet\varphi+\frac{1}{2} EX\bullet(T\bullet\varphi+f\varphi) \nonumber\\&&-\frac{i}{2} \tau X\bullet\varphi-\frac 12 f\tau JX\bullet\varphi.
\label{c2}
\end{eqnarray}
Take Equation (\ref{c2}) and substract (\ref{c1}) to get
$$-X(f)\varphi + fEX\bullet\varphi -g(T, EX) \varphi - \nabla_XT\bullet\varphi -\frac{\tau}{2} T\bullet JX\bullet\varphi =0.$$
Taking the real part of the scalar product of the last equation with $\varphi$ and using that $<iX\bullet\varphi, \varphi> = -g(T, JX) \vert\varphi\vert^2$, we get
$$X(f) = - g(T, EX) + \tau g(JX, T).$$
The imaginary part of the same scalar product gives $\nabla_XT= f(EX-\tau JX),$ which gives that there exists an immersion $F$ from $M$ into $\Ekt$ with shape operator $dF\circ E \circ dF^{-1}$ and $\xi = dF(T) + f\nu$.\\\\
Now, we state the main result of this section, which characterize any isometric immersion of a surface $(M,g)$ into 
$\Ekt$.\begin{thm}\label{thmEkt}
Let $\kappa,\tau\in\R$ with $\kappa-4\tau^2\neq0$. Consider $(M^2, g)$  a simply connected oriented Riemannian surface. We denote by $E$ a field of symmetric endomorphisms of $TM$, with trace equal to $2H$. The following statements are equivalent:
\begin{enumerate}
\item There exists an isometric immersion $F$ of $(M^2,g)$ into $\Ekt$ with shape operator $E$, mean curvature $H$ and such that, over $M$, the vertical vector is $\xi=dF(T)+f\nu$, where $\nu$ is the unit normal vector to the surface, $f$ is a real function on $M$ and $T$ the tangential part of $\xi$.
\item There exists a $\spinc$ structure on $M$ carrying a non-trivial spinor field $\varphi$ satisfying
$$\nabla_X\varphi=-\frac{1}{2}EX\bullet\varphi+i\frac{\tau}{2}X\bullet\overline{\varphi}.$$
Moreover, the auxiliary bundle has a connection of curvature given, in any local orthonormal frame $\{t_1, t_2\}$, by   $\Omega(t_1, t_2)= -(\kappa-4\tau^2)f =-(\kappa-4\tau^2) \frac{<\varphi, \overline \varphi>}{\vert\varphi\vert^2} $.
\item There exists a $\spinc$ structure on $M$ carrying a non-trivial spinor field $\varphi$ of constant norm satisfying
$$D\varphi=H\varphi-i\tau\overline{\varphi}.$$
Moreover, the auxiliary bundle has a connection of curvature given, in any local orthonormal frame $\{t_1, t_2\}$, by  $\Omega(t_1, t_2)= -(\kappa-4\tau^2)f =-(\kappa-4\tau^2) \frac{<\varphi, \overline \varphi>}{\vert\varphi\vert^2} $.
\end{enumerate} 
\end{thm}
{\bf Proof:} Proposition \ref{propEkt} and Lemma \ref{lemEkt} give the equivalence between the  first two statements. If the statement ($2$) holds, it is easy to check that in this case the Dirac operator acts on $\varphi$ to give 
$D\varphi = H\varphi - i\tau \overline\varphi.$ Moreover, for any $X\in \Gamma(TM)$, we have
\begin{eqnarray*}
 X(\vert\varphi\vert^2) &=& 2\Re\<\nabla_X \varphi, \varphi\> \\&=&
\Re \<i\tau X\bullet\overline \varphi, \varphi\> =0.
 \end{eqnarray*}
Hence $\varphi$ is of constant norm . Now, consider a non-trivial
spinor field $\varphi$ of constant length, which satisfies $D \varphi  =H \varphi -i \tau \overline \varphi.$ Define the following $2$-tensors on $(M^2,g)$
$$T^\varphi_{\pm}(X,Y)=\Re\<\nabla_X\varphi^\pm,Y\bullet\varphi^\mp\>\;.$$ First note that 
\begin{equation}
  \label{eq:trE}
 \mathrm{tr}T^\varphi_{\pm}=-\Re\<D\varphi^\pm,\varphi^\mp\>= -H|\varphi^\mp|^2\;. 
\end{equation}
Moreover, we have the following relations \cite{Mo}
\begin{eqnarray}
 T^\varphi_{\pm}(t_1,t_2) &=& \tau |\varphi^\mp|^2+ T^\varphi_{\pm}(t_2,t_1),
\label{111}
\end{eqnarray}
\begin{eqnarray}
\nabla_X\varphi^+ &=& \frac{T^\varphi_+(X)}{|\varphi^-|^2}\bullet\varphi^-,\\
\gn_X\varphi^-&=&\frac{T^\varphi_-(X)}{|\varphi^+|^2}\bullet\varphi^+\;,\\
|\varphi^+|^2T^\varphi_+&=&|\varphi^-|^2T^\varphi_-,
\end{eqnarray}
where the vector field $T^\varphi_+(X)$ is defined by
$g(T^\varphi_+(X),Y)=T^\varphi_+(X,Y)$ for  $Y\in \Gamma(TM).$ Now let $F^\varphi:=T^\varphi_+ + T^\varphi_-$. Thus, we have 

$$\frac{F^\varphi}{|\varphi|^2}=\frac{T^\varphi_+}{|\varphi^\varphi_-|^2}=\frac{T^\varphi_-}{|\varphi^+|^2}\;.$$ Hence
$F^\varphi/|\varphi|^2$ is well defined on the whole surface $M$, and \begin{equation}\label{preqrsk}\nabla_X\varphi=\gn_X\varphi^++\gn_X\varphi^-=\frac{F^\varphi(X)}{|\varphi|^2}\bullet\varphi,\end{equation}
 where the vector field $F^\varphi(X)$ is defined by
$g(F^\varphi(X),Y)=F^\varphi(X,Y)$, for all $Y\in \Gamma(TM).$ Note that by Equation
(\ref{111}), the $2$-tensor $F^\varphi$ is not symmetric. Define now the symmetric $2$-tensor
$$T^\varphi(X,Y)=-\frac{1}{2|\varphi|^2}\left(F^\varphi(X,Y)+F^\varphi(Y,X)\right)\;.$$ It is straigthforward to
show that
$$T^\varphi(t_1,t_1)=-F^\varphi(t_1,t_1)/|\varphi|^2\quad,\quad T^\varphi(t_2,t_2)=-F^\varphi(t_2,t_2)/|\varphi|^2\;,$$
$$T^\varphi(t_1,t_2)=-F^\varphi(t_1,t_2)/|\varphi|^2+\frac {\tau}{2}\quad\text{and}\quad T^\varphi(t_2,t_1)=-F^\varphi(t_2,t_1)/|\varphi|^2-\frac {\tau}{2}.$$ 
Taking into account these last relations in Equation
(\ref{preqrsk}), we conclude 
\begin{eqnarray*}
 \nabla_X\varphi &=& -T^\varphi(X)\bullet\varphi + i \frac{\tau}{2} X\bullet\overline \varphi.
\end{eqnarray*}
\subsection{Application: a spinorial proof of Daniel correspondence}
In \cite{Da2}, B. Daniel gave a  Lawson type correspondence for constant mean curvature surfaces in $\Ekt$. Namely, he proved the following
\begin{thm}\label{thmdan}
Let $\Ekto$ and $\Ektt$ be two $3$-dimensional homogeneous manifolds  with four dimensional isometry group and assume that $\kappa_1-4\tau_1^2=\kappa_2-4\tau_2^2$. Consider $\xi_1$ and $\xi_2$ be the vertical vectors of $\Ekto$ and $\Ektt$ respectively and $(M^2,g)$  a simply connected surface isometrically immersed into $\Ekto$ with constant mean curvature $H_1$ so that $H_1^2\geq\tau_2^2-\tau_1^2$. We denote by  $\nu_1$ be the unit inner normal of the immersion, $T_1$ the tangential projection of $\xi_1$ and $f=<\nu_1,\xi_1>$. \\
Let $H_2\in \R$ and $\theta\in\R$ so that 
$$H_2^2+\tau_2^2=H_1^2 + \tau_1^2,\ \ \ \ \ \text{and}\ \ \ \ \ \tau_2+iH_2=e^{i\theta}(\tau_1+iH_1).$$
Then, there exists an isometric immersion $F$ from $(M^2,g)$ into $\Ektt$ with mean curvature $H_2$ and so that over $M$
$$\xi_2=dF(T_2)+f\nu_2,$$
where $\nu_2$ is the unit normal inner vector of the immersion and $T_2$ the tangential part of $\xi_2$. Moreover, the respective shape operator $E_1$ and $E_2$ are related by the following
$$E_2-H_2\iid=e^{\theta J}(E_1-H_1\iid).$$
\end{thm}

With the help of Theorem \ref{thmEkt}, we give an alternative proof of this results using spinors.\\ \\
{\bf Proof of Theorem \ref{thmdan}:} Since $M^2$ is isometrically immersed into $\Ekto$ there exists a spinor field $\varphi_1$ of constant norm (say $|\varphi_1|=1$) satisfying
$$D\varphi_1=H_1\varphi_1-i\tau\overline{\varphi_1},$$
associated with the $\spinc$ structure whose line bundle has a connection of curvature given by  $\Omega= -(\kappa-4\tau^2)f $, where $f= \frac{<\varphi, \overline\varphi>}{\vert\varphi\vert}$. We deduce that 
$$D\varphi_1^+=H_1\varphi_1^-+i\tau_1\varphi_1^-\,$$
$$D\varphi_1^-=H_1\varphi_1^+-i\tau_1\varphi_1^+.$$
Now, we define $\varphi_2=\varphi_1^++e^{i\theta}\varphi_1^-$.
First, we have
\begin{eqnarray*}
 D\varphi_2&=&D\varphi_1^++e^{i\theta}D\varphi_1^-\\
&=&(H_1+i\tau_1)\varphi_1^--ie^{i\theta}(\tau_1+iH_1)\varphi^+_1
\end{eqnarray*}
Since $\tau_2+iH_2=e^{i\theta}(\tau_1+iH_1)$, we deduce that $H_1+i\tau_1=e^{i\theta}(H_2+i\tau_2)$ and so $D\varphi_2=H_2\varphi_2-i\tau_2\overline{\varphi_2}$. Secondly, 
$$\frac{<\varphi_1, \ov\varphi_1>}{\vert\varphi_1\vert^2} =  \frac{<\varphi_2, \ov\varphi_2>}{\vert\varphi_2\vert^2}.$$
Now, since $\kappa_1-4\tau_1^2=\kappa_2-4\tau_2^2$, the considered $\spinc$ structure on $M$ is given by $i\Omega=-i(\kappa_2-4\tau_2^2)f$ and hence, by Theorem \ref{thmEkt}, there exists an isometric immersion $F$ from $(M^2,g)$ into $\Ektt$ with mean curvature $H_2$ and so that $\xi_2=dF(T_2)+f \nu_2$, where $\nu_2$ is the unit normal inner vector of the surface and $T_2$ the tangential part of $\xi_2$.
\begin{rem}
By the proof of Proposition \ref{propEkt},  we have that
$$g(T_2,t_1)|\varphi_2|^2=\Re\<it_2\bullet\varphi_2,\varphi_2\> \quad\text{and}\quad g(T_2,t_2)|\varphi_2|^2=-\Re\<it_1\bullet\varphi_2,\varphi_2\>.$$
So, it is easy to see that $T_2 =e^{\theta J}(T_1)$.
\end{rem}
\section{Isometric immersions into $\CM$ via spinors}
In this section, we consider the canonical $\spinc$ structure on $\CM$ carrying a parallel spinor field $\psi$ lying in $\Sigma^+ (\CM)$. The restriction of this $\spinc$ structure to any hypersurface $M^3$ defines a $\spinc$ structure on $M$ with a special spinor field. This spinor field characterizes the isometric immersion of $M$ into $\CM$.
\subsection{Special spinors fields on $\CM$ and their surfaces}
Assume that there exists an isometric immersion of $(M^3,g)$ into $\CM$ with shape operator $II$. By section \ref{cm}, we know that $M$ has an almost contact metric structure $(\Chi,\xi,\eta)$ such that $\Chi X = JX -\eta(X) \nu$ for every $X\in \Gamma(TM)$.
\begin{lem}
 The restriction $\varphi$  of the parallel spinor $\psi$ on $\CM$ is a solution of the generalized Killing equation
\begin{eqnarray}
 \nabla_X\varphi + \frac 12 II X\bullet\varphi=0,
\end{eqnarray}
Moreover, $\varphi$ satisfies $\xi\bullet\varphi = -i\varphi$. The curvature 2-form of the auxiliary line bundle associated with the induced $\spinc$ structure is given by $\Omega (X, Y) = -6c \ltimes (X, Y)$, where $\ltimes$ is the K\"{a}hler form of $\CM$ given by $\ltimes(X, Y) = \ov g(JX, Y)$.
\label{lema2}
\end{lem}
{\bf Proof:} First, since $\psi$ is parallel, we have $ D^{\CM} \psi = (\nabla^{\CM})^* \nabla^{\CM }\psi =0$. Hence, by the Schr\"{o}dinger-Lichnerowicz formula, we get 
\begin{eqnarray}
\Omega^{\CM} \cdot\psi = 12 c i \psi.
\end{eqnarray}
By the Gauss formula (\ref{spingauss}), the restriction $\varphi$ of the parallel spinor $\psi$ on $\CM$ satisfies
$$\nabla_X\varphi =- \frac 12 II X\bullet\varphi.$$
Since the spinor $\psi$ is parallel, Equality (\ref{RRicci}) gives
$$\mathrm{Ric}^{\CM}(X)\cdot\psi = i (X\lrcorner \Omega^{\Ekt})\cdot\psi $$
Where $\Ric$ is the Ricci tensor of $\CM$. Therfore, we compute,
\begin{eqnarray*}
(\nu\lrcorner \Omega^{\CM})\bullet\varphi &=&(\nu\lrcorner \Omega^{\CM})\cdot\nu\cdot\psi|_M \\ &=&
-\nu\cdot(\nu\lrcorner \Omega^{\CM})\cdot\psi|_M\\&=& i  \nu \cdot\mathrm{Ric}^{\CM}\ \nu\cdot\psi_{\vert_M}\\&=& -6ci\varphi.
\end{eqnarray*}
By Equation (\ref{glucose}), we get that 
\begin{eqnarray}
 \Omega\bullet \varphi = 6ci\varphi.
\label{xi}
\end{eqnarray}
Now, for any $X, Y \in \Gamma(TM)$, we have 
$$\Omega (X, Y) = \Omega^{\CM} (X, Y)= -\rho(X, Y)=-\mathrm{Ric} (JX, Y) = -6c \ov g(JX, Y).$$
Let $e_1$ be a unit vector field tangent to $M$ such that $\{e_1, e_2 =Je_1, \xi\}$ is an orthonormal basis of $TM$. In this basis, we have 
$$\Omega\bullet\varphi =\Omega(e_1, e_2)\  e_1\bullet e_2 \bullet \varphi+ \Omega(e_1, \xi)\  e_1\bullet \xi \bullet \varphi +\Omega(e_2, \xi)\  e_2\bullet \xi \bullet \varphi.$$
But, 
$$\Omega(e_1, e_2)=  -6c\ \ \ \ \ \ \ \ \text{and}\ \ \ \ \ \ \Omega(e_1, \xi) = \Omega(e_2, \xi) = 0.$$
Finally, $\Omega\bullet\varphi = -6c e_1 \bullet e_2 \bullet\varphi$. Using (\ref{xi}) and the fact that $e_1 \bullet e_2 \bullet \xi \bullet \varphi = - \varphi$, we conclude that $\xi \bullet \varphi = -i \varphi$.
\begin{lem}
Let $E$ be a field of symmetric endomorphisms on a $\spinc$ manifold $M^3$ of dimension $3$, then 
\begin{eqnarray}
 E(e_i)\bullet E(e_j) - E(e_j)\bullet E(e_i )&=& 2 (a_{j3} a_{i2} - a_{j2}a_{i3}) e_1 \nonumber \\&& +2(a_{i3}a_{j1} - a_{i1}a_{j3}) e_2 \nonumber \\ && + 2(a_{i1}a_{j2} - a_{i2}a_{j1})e_3,
\end{eqnarray}
where $(a_{ij})_{i,j}$ is the matrix of $E$ written in any local orthonormal frame of $TM$.
\label{aij}
\end{lem}
\begin{prop}
Let $(M^3,g)$ be a Riemannian $\spinc$ manifold endowed with an almost contact metric structure $(\Chi,\xi,\eta)$. Assume that there exists a non-trivial spinor $\varphi$  satisfying
$$\nabla_X\varphi=-\frac{1}{2}EX\bullet\varphi \ \ \ \text{and}\ \ \ \ \xi\bullet\varphi=-i\varphi,$$
where $E$ is a field of symmetric endomorphisms on $M$. We suppose that the curvature 2-form of the connection 
on the auxiliary line bundle associated with the $\spinc$ structure is given by $\Omega(e_1,e_2)=-6c$ and
 $\Omega(e_i,e_j)=0$ elsewhere in the basis $\{e_1,e_2=\Chi e_1,e_3=\xi\}$.
Hence,  the Gauss equation for $\CM$  is satisfied if and only if the Codazzi equation for $\CM$ is satisfied.
\label{GGCC}
\end{prop}
{\bf Proof:} We compute the spinorial curvature $\mathcal{R}$ on $\varphi$, we get
$$\mathcal{R}_{X, Y} \varphi= -\frac 12 d{^\nabla} E(X, Y)\bullet\varphi + \frac 14 (EY\bullet EX - EX \bullet EY)\bullet\varphi.$$
In the basis $\{e_1, e_2 = \Chi e_1, e_3 = \xi\}$, the Ricci identity (\ref{RRicci}) gives that 
\begin{eqnarray*}
 \frac 12 \mathrm{Ric}(X)\bullet\varphi -\frac i2 (X\lrcorner\Omega)\bullet\varphi &=& \frac 14 \sum_{k=1}^3 e_k\bullet(EX\bullet Ee_k - Ee_k\bullet EX)\bullet\varphi \\&&
-\frac 12 \sum_{k=1}^3  e_k\bullet d^{\nabla} E(e_k, X )\bullet\varphi.
\end{eqnarray*}
By Lemma \ref{aij} and for $X= e_1$, the last identity  becomes   
\begin{eqnarray}
 &\ &(\mathrm{R}_{1221} +\mathrm{R}_{1331} -a_{11}a_{33} - a_{11}a_{22} + a_{13}^2 +a_{12}^2- 5c )e_1 \bullet \varphi \nonumber\\ && +( \mathrm{R}_{1332} -a_{12}a_{33} + a_{32}a_{13})e_2\bullet \varphi \nonumber \\ && + (\mathrm{R}_{1223} -a_{22}a_{13} + a_{32}a_{12})e_3\bullet \varphi \nonumber\\ &=&
- e_2\bullet d^{\nabla}E (e_2, e_1)\bullet\varphi- e_3\bullet d^{\nabla}E (e_3, e_1)\bullet\varphi  \nonumber\\ && +c e_1\bullet\varphi.
\label {e1varphi}
\end{eqnarray}
Since $\vert \varphi\vert$ is constant ($\vert \varphi \vert =1$), the set $\{\varphi, e_1\bullet\varphi, e_2 \bullet\varphi, e_3\bullet\varphi\}$ is an orthonormal frame of $\Sigma M$ with respect to the real scalar product $\pre \< .,. \>$. Hence, from Equation (\ref{e1varphi}) we deduce
\begin{eqnarray*}
 \mathrm{R}_{1221} +\mathrm{R}_{1331} -(a_{11}a_{33} + a_{11}a_{22} - a_{13}^2 -a_{12}^2+ 5c) & =& g(d^\nabla E(e_1, e_3), e_3) - g(d^\nabla E(e_1, e_3), e_2) + c\\ 
\mathrm{R}_{1332} -(a_{12}a_{33} - a_{32}a_{13}) &=& g(d^\nabla E(e_1, e_3), e_1)\\
\mathrm{R}_{1223} -(a_{22}a_{13} - a_{32}a_{12})&=& g(d^\nabla E(e_1, e_2), e_1)\\
 g(d^\nabla E(e_1, e_2), e_2) &=& - g(d^\nabla E(e_1, e_3), e_3) 
\end{eqnarray*}
The same computation holds for the unit vector fields $e_2$ and $e_3$ and we get
\begin{eqnarray*}
\mathrm{R}_{2331}  -(a_{12}a_{33} - a_{13}a_{23})& =& - g(d^\nabla E(e_2, e_3), e_2)\\
\mathrm{R}_{2332} + \mathrm{R}_{2112}  -(a_{22}a_{33} + a_{22}a_{11} - a_{13}^2 - a_{12}^2 +5c) &=& g(d^\nabla E(e_2, e_3), e_1) +g(d^\nabla E(e_1, e_2), e_3) +c\\
\mathrm{R}_{2113} -(a_{23}a_{11} - a_{12}a_{13}) &=& - g(d^\nabla E(e_1, e_2), e_2)\\
 g(d^\nabla E(e_1, e_2), e_1) &=&  g(d^\nabla E(e_2, e_3), e_3)\\
\mathrm{R}_{3221}  -(a_{13}a_{22} - a_{23}a_{21})  &=& - g(d^\nabla E(e_2, e_3), e_3)\\
\mathrm{R}_{3112}  - (a_{32}a_{11} - a_{31}a_{12}) &=& g(d^\nabla E(e_1, e_3), e_3)\\
\mathrm{R}_{3113}+ \mathrm{R}_{3223}-(a_{22}a_{33} - a_{11}a_{33} + a_{13}^2 + a_{23}^2) &= & g(d^\nabla E(e_2, e_3), e_1) - g(d^\nabla E(e_1, e_3), e_2)\\
g(d^\nabla E(e_2, e_3), e_2) &=& - g(d^\nabla E(e_1, e_3), e_1)
\end{eqnarray*}
The last twelve equations imply that  the Gauss equation for $\CM$ is satisfied if and only if the Codazzi equation
 for $\CM$ is satisfied.
\subsection{Spinorial characterization of hypersurfaces of $\CM$}
Now, we give the main result of this section:
\begin{thm}\label{thmCM}
Let $(M^3,g)$ be a simply connected oriented Riemannian manifold endowed with an almost contact metric structure
 $(\Chi,\xi,\eta)$. Let $E$ be a field of symmetric endomorphisms on $M$ with trace equal to $3H$. Assume that the Gauss or the 
Codazzi equation for $\CM$ is satisfied. Then, the following statements are equivalent:
\begin{enumerate}
\item There exists an isometric immersion of $(M^3,g)$ into $\CM$ with shape operator $E$, mean curvature $H$ and so that, over $M$, the complex structure of $\CM$ is given by $J=\Chi+\eta(\cdot)\nu$, where $\nu$ is the unit normal vector of the immersion.
\item There exists a $\spinc$ structure on $M$ carrying a  non-trivial spinor $\varphi$  satisfying
$$\nabla_X\varphi=-\frac{1}{2}EX\bullet\varphi\ \ \ \text{and}\ \ \  \xi\bullet\varphi= - i\varphi.$$
The curvature 2-form of the connection on  the auxiliary bundle associated with the $\spinc$ structure is given by $\Omega(e_1,e_2)=-6c$ and $\Omega(e_i,e_j)=0$ elsewhere in the basis $\{e_1,e_2=\Chi e_1,e_3=\xi\}$.
\item There exists a $\spinc$ structure on $M$ carrying a  non-trivial spinor $\varphi$ of constant norm and satisfying
$$D\varphi=\frac{3}{2}H\varphi\ \ \ \text{and} \ \ \ \xi\bullet\varphi= - i\varphi.$$
The curvature 2-form of the connection on the auxiliary bundle associated with the $\spinc$ structure is given by $\Omega(e_1,e_2)=-6c$ and $\Omega(e_i,e_j)=0$ elsewhere in the basis $\{e_1,e_2=\Chi e_1,e_3=\xi\}$.
\end{enumerate}
\end{thm}
{\bf Proof:} By Lemma \ref{lema2}, the first statement implies the second one. Using Proposition \ref{GGCC}, to show that $2 \Longrightarrow 1$, it suffies to show that $\nabla_X \xi = \Chi EX$. In fact, we simply compute the derivative of $\xi\bullet\varphi = -i \varphi$ in the direction of $X\in \Gamma (TM)$ to get
\begin{eqnarray*}
 \nabla_X \xi \bullet\varphi &=& \frac i2 EX\bullet\varphi + \frac 12 \xi\bullet EX\bullet\varphi 
\end{eqnarray*}
Using that $-i e_2 \bullet\varphi = e_1 \bullet\varphi$, the last equation reduces to 
$$\nabla_X\xi \bullet\varphi - g(EX, e_1) e_2\bullet\varphi  + g(EX, e_2) e_1\bullet\varphi =0.$$
Finally $\nabla_X \xi\bullet\varphi = \Chi EX$. Now, we compute the derivative of $-i e_2\bullet\phi = e_1\bullet\phi$ in the direction of $e_1$ to get
$$\nabla_{e_1} (\Chi e_1)\bullet\phi -\frac 12 e_2\bullet Ee_1\bullet \phi = i \nabla_{e_1}e_1 \bullet \phi -\frac i2 e_1\bullet Ee_1\bullet\phi.$$
But, using that $\xi\bullet\phi = -i\phi$, we have
$$\frac 12 e_2\bullet Ee_1\bullet \phi -\frac i2 e_1\bullet Ee_1\bullet\phi = -a_{11} \xi\bullet\phi -a_{12} \phi.$$
Denoting by $\Gamma_{ij}^k$ the Christoffel symbols of $\{e_1, \Chi e_1, \xi\}$, we have $\nabla_{e_1}e_1 = \Gamma_{11}^1 e_1 + \Gamma_{11}^2 e_2 + \Gamma_{11}^3 e_3$. Moreover, using that $\nabla_{e_1}e_3 = \Chi Ee_1$, we get
$$\Gamma_{11}^3 = g(\nabla_{e_1}e_1, e_3) = -g (e_1, \nabla_{e_1}e_3) = a_{12}.$$
Hence, $\nabla_{e_1} (\Chi e_1)\bullet\phi = -a_{11} \xi\bullet\phi +  \Gamma_{11}^1 e_2\bullet\phi +  \Gamma_{11}^2 e_2\bullet\phi.$ Finally
$$\nabla_{e_1} (\Chi e_1)\bullet\phi - \Chi(\nabla_{e_1}e_1)\bullet\phi = -a_{11} \xi\bullet\phi,$$
which is Equation (\ref{cond3}) for $X = Y =e_1$. Similary, 
we compute the derivative of $-i e_2\bullet\phi = e_1\bullet\phi$ in the direction of $e_2$ and $\xi$ to get Equation 
(\ref{cond3}) for any $X, Y\in \Gamma(TM)$. It is easy to see that the assertion 2 implies the assertion 3. 
For $3 \Rightarrow 2$, since $\varphi$ is of constant norm ($\vert \varphi \vert =1$), the set 
$\{\varphi, e_1\bullet\varphi,  e_2\bullet\varphi,  e_3\bullet\varphi\}$ is a local orthonormal frame of
 $\Sigma M$ with respect to the real scalar product $\pre \<., .\>$. Hence, for every $X\in \Gamma(TM)$, we have
\begin{eqnarray} 
 \nabla_X\varphi = \eta(X) \varphi + \ell(X)\bullet\varphi,
\label{eq1}
\end{eqnarray}
where $\eta$ is a $1$-form and $\ell$ is a $(1,1)$-tensor field. Moreover it is easy to check that
 $\eta =\frac{d(\vert\varphi\vert^2)}{2\vert\varphi\vert^2}$ and $\ell(X)=- \ell^\varphi(X)$. Since $\varphi$ is 
of constant norm we have $\eta=0$. Moreover, $\ell(X)=- \ell^\varphi(X)$  is symmetric of trace $\frac 32 H$. 
It suffices to consider $E = 2\ell^\varphi$ to get the second assertion.
\subsection{Characterization of  Sasaki hypersurfaces}
Theorem \ref{thmCM}  characterizes isometric immersions of almost contact metric manifolds into $\CM$ providing that the shape operator $E$ satisfies the Gauss or the Codazzi equation for $\CM$. 
In this subsection, we eliminate this restriction and we replace it by some geometric conditions on the almost contact metric manifold.\\\\
In Section $3$, we showed that the $3$-dimensional homogeneous manifolds $\Ekt$ ($\tau \neq 0, \kappa-4\tau^2 \neq 0$), which are Sasaki, have a $\spinc$ structure 
(the canonical $\spinc$ structure) carrying a Killing spinor field $\varphi$ of Killing constant $\frac{\tau}{2}$. Moreover $\xi\cdot\varphi = -i\varphi$ and 
\begin{eqnarray}
\Omega^{\Ekt} (e_1, e_2) = -(\kappa - 4 \tau^2) \ \ \text{and} \ \ \ \Omega^{\Ekt}(e_i, e_j) = 0,
\end{eqnarray}
in the basis $\{e_1, \Chi e_1 = e_2, e_3 = \xi\}$. Hence, the statement ($2$) of Theorem \ref{thmCM} is satisfied for 
$E = -\tau\ \mathrm{Id}$ and $c= \frac{\kappa-4\tau^2}{6} \neq 0$. But $\Ekt$ cannot be immersed into $\CM$ 
($c =\frac{\kappa-4\tau^2}{6} \neq 0$) with second fundamental form $E = -\tau \ \mathrm{Id}$ because we know that totally umbilic hypersurfaces in $\CM$ cannot exist. Moreover, the Codazzi equation is not satisfied. In fact, it is easy to check that $d^\nabla E(e_1, e_2) = 0,$ and
$$ c\{ \eta(e_1)\Chi e_2 - \eta(e_2)\Chi e_1 + 2 g(e_1, \Chi e_2)\xi\} = -2 c \xi \neq 0 .$$ 
From this example, it is clear that the condition ``$E$ satisfies the Gauss equation or the Codazzi equation'' is a 
necessary condition to immerse in $\CM$ an almost contact metric manifold $M$ satisfying the statement ($2$) of
Theorem \ref{thmCM} and even if the manifold $M$ is Sasaki. However, we can state the following:
\begin{thm}\label{sasa}
 Let $(M^3,g)$ be a simply connected oriented Riemannian manifold endowed with a Sasakian structure
 $(\Chi,\xi,\eta)$. Then, the following statements are equivalent:
\begin{enumerate}
\item There exists an isometric immersion of $(M^3,g)$ into $\CM$ with mean curvature $H$ and so that, over $M$, the complex structure of $\CM$ is given by $J=\Chi+\eta(\cdot)\nu$, where $\nu$ is the unit normal vector of the immersion. 
\item There exists a $\spinc$ structure on $M$ carrying a  non-trivial spinor $\varphi$  satisfying
$$\nabla_X\varphi=-\frac{1}{2}X\bullet\varphi -\frac{i}{2} c\ \eta(X)\varphi\ \ \ \text{and}\ \ \  \xi\bullet\varphi= - i\varphi.$$
The curvature 2-form of the connection on  the auxiliary bundle associated with the 
$\spinc$ structure is given by $\Omega(e_1,e_2)=-6c$ and $\Omega(e_i,e_j)=0$ elsewhere 
in the basis $\{e_1,e_2=\Chi e_1,e_3=\xi\}$.
\end{enumerate}
In this case, $M$ is of constant mean curvature $H = \frac{3-c}{3}$ and the shape operator $E$ is given by $E = \iid -c \eta(\cdot)\xi$.
\end{thm}
{\bf Proof:} Assume that $(M^3,g)$ is a Sasaki manifold immersed into $\CM$ with shape operator $E$. Since
 $\xi$ is a Killing vector field, Equation (\ref{condd4}) implies that $\Chi (E\xi ) = \nabla_\xi \xi =0$ and hence 
$E\xi = f\xi$, where $f$ is a real function on $M$. Also, from Equation (\ref{condd4}) and since $\nabla_X\xi = \Chi X$, we get
$\Chi(EX-X) =0$, for all $X\in \Gamma(TM)$. Then, $$EX- X = g(EX-X, \xi) \xi.$$
But, $g(EX-X, \xi) = (f-1) g(X, \xi)$ which gives that  $EX = X + (f-1) g(X, \xi)\xi.$ It is straightforward to check that
\begin{eqnarray*}
(\nabla_X E)(Y) - (\nabla_Y E)(X) &=& -(f - 1) \{ \eta (X) \Chi Y-  \eta(Y)\Chi X + 2g(X, \Chi Y) \xi \}\\
&& + \{ df(X) \eta(Y) - df(Y)\eta(X) \}\xi,
\end{eqnarray*}
for all vectors $X, Y \in \Gamma(TM)$. Comparing the last equation with (\ref{codazziCMM}), we get $f-1 = -c$. This
 gives $EX = X -c\ \eta(X)\xi$ and by Theorem   \ref{thmCM}, we get the statement ($2$). Now, we assume that the statement ($2$) holds, i.e., we have 
on $M$  a $\spinc$ structure  carrying a  non-trivial spinor $\varphi$ satisfying
\begin{eqnarray}\label{gene-killing}
 \nabla_X\varphi=-\frac{1}{2}X\bullet\varphi -\frac{i}{2} c\ \eta(X)\varphi\ \ \ \text{and}\ \ \  \xi\bullet\varphi= - i\varphi.
\end{eqnarray}
The curvature 2-form of the connection on  the auxiliary bundle associated with the $\spinc$ structure is given by $\Omega(e_1,e_2)=-6c$ and $\Omega(e_i,e_j)=0$ elsewhere in the basis $\{e_1,e_2=\Chi e_1,e_3=\xi\}$.  
We denote by $E$ the endomorphism given for all $X\in \Gamma(TM)$, by  $EX = X-c\eta(X)\xi $. From (\ref{gene-killing}), 
we have $\nabla_X\varphi=-\frac{1}{2}EX\bullet\varphi$ and  we can check that $E = \mathrm{Id} -c\eta(\cdot)\xi$ satisfies, for all vectors $X, Y \in \Gamma(TM)$, 
$$(\nabla_X E)(Y) - (\nabla_Y E)(X) = c \{ \eta (X) \Chi Y-  \eta(Y)\Chi X + 2g(X, \Chi Y) \xi \},$$
which is the Codazzi equation (\ref{codazziCM}). By Theorem \ref{thmCM}, $M$ is immersed
 into $\CM$ with shape operator $E$. Additionally, since $EX = X-c\eta(X)\xi $, we have $H =\frac{3-c}{3}$.
\begin{rem}
From the above example, $\Ekt$ with $\tau \neq 0$ endowed with their canonical $\spinc$ structure cannot 
be immersed into $\CM$ for $c= \frac{\kappa-4\tau^2}{6} \neq 0$. In fact, the Killing spinor of Killing constant 
$\frac{\tau}{2}$
 does not satisfy assertion ($2$) of Theorem \ref{sasa} because for example, when $\tau = -1$,  the endomorphism 
$E = \mathrm{Id}$ 
is not of the form $E  = \iid -c\ \eta(\cdot)\xi$. On the other side, it is known that there exists an isometric embedding of $\Ekt$, $\tau\neq0$, into 
$\M^2(\frac {\kappa}{4} - \tau^2)$ of constant mean curvature $H = \frac{\kappa-16\tau^{2}}{12\tau}$ \cite{CMC}. 
In a recent work \cite{nakadroth}, the authors used the canonical and the anti-canonical $\spinc$ structures on $\Ekt$, 
to define another $\spinc$ structure on $\Ekt$  satisfying assertion ($2$) of Theorem \ref{sasa} and hence allowing to immerse $\Ekt$ into $\CM$. Other geometric applications are also given.
\end{rem}
{\bf Acknowledgement:} Both authors are grateful to Oussama Hijazi for his encouragements, valuable comments and relevant remarks.

\end{document}